\documentclass[12pt,twoside]{article} 
\usepackage[cp1251]{inputenc}
\usepackage{latexsym,amsfonts,amssymb,amsmath}
\usepackage{euscript}
\usepackage[T2A]{fontenc}
\usepackage{amsthm}
\usepackage[english]{babel}
\usepackage{mathrsfs}

\begin{document}


\newcommand{\ang}[1]{\langle{#1}\rangle}

\thispagestyle{empty}

\noindent {UDC 519.41/47}

\setcounter{page}{1}

{\small\vspace{1mm}\noindent{\bf N.S.~Chernikov} (In-te Math. of NAS
of Ukraine, Kyiv, Ukraine)

\vspace{1mm}\noindent{\bf A NOTE ON GROUPS WITH THE MINIMAL
CONDI\-TIONS FOR NONABELIAN AND ABELIAN SUBGROUPS}

\noindent chern@imath.kiev.ua

\vspace{0.5cm}\noindent{\small We give a new proof of the known
Shunkov's Theorem on locally finite groups with the minimal
condition for nonabelian subgroups and also an extension of the
known Suchkova-Shunkov Theorem on Shunkov groups with the minimal
condition for abelian subgroups.}

\thispagestyle{empty}
\vspace{8mm}

\noindent The celebrated Shunkov's Theorem \cite{1} asserts that a
locally finite group, which satisfies the minimal condition for
abelian subgroups, is Chernikov. Another known Shunkov's Theorem
\cite{2} establishes that a nonabelian locally finite group $H$,
satisfying the minimal condition for nonabelian subgroups, is
Chernikov too. Further, a nonabelian almost solvable group with the
minimal condition for nonabelian subgroups is itself Chernikov
(S.N.Chernikov's Theorem \cite{3}).

Below we'll give a new proof of Shunkov's Theorem \cite{2}, which is
based on Shunkov's Theorem \cite{1} and S.N.Chernikov's Theorem
\cite{3} and on the following proposition.

\vskip2mm \textbf{Proposition.} \textit{If an abelian subgroup $A$
of the locally finite group $G$ is nonnormal and maximal in it, then
\begin{equation}\label{1}
A\cap  A^{g}=Z(G),\quad \forall g\in G\backslash A,
\end{equation}and for some normal subgroup $N$ of $G$,
\begin{equation}\label{2}
    G=AN \mbox{ and } A\cap N=Z(G).
\end{equation}}

{\itshape\bfseries Proof. }  Let $A$ be nonnormal and maximal in
$G$. It is easy to see that for $g\in G\backslash A$,
$G=\ang{A,A^{g}}$. So $A\cap A^{g}\subseteq Z(G)$. But, clearly,
$Z(G)\subseteq A,A^{g}$. Thus \eqref{1} is correct. Therefore
\begin{equation}\label{3}
    A/Z(G) \cap (A/Z(G))^{g}=1,\quad \forall g\in G/Z(G)\backslash
    A/Z(G).
\end{equation}Then, with regard to \eqref{3}, by
Busarkin-Starostin-Kegel Theorem \cite{4}-\cite{6}, for some normal
subgroup $N$ of $G$, \eqref{2} are correct.

\vskip2mm {\itshape\bfseries Proof of Shunkov's Theorem} \cite{2}.
Let $H$ be non-Chernikov. Then, obviously, $H$ contains some non-
(Chernikov or abelian) subgroup $G$ such that any its proper
subgroup is Chernikov or abelian. In view of Shunkov's Theorem
\cite{1}, some maximal abelian subgroup $A$ of $G$ is non-Chernikov.
Obviously, $A$ is maximal in $G$. In view of S.N.Chernikov's Theorem
\cite{3}, $G$ is not almost solvable. Hence follows that $A$ is not
normal in $G$. Consequently, for some proper normal subgroup $N$ of
$G$, $G=AN$ (see Proposition). Clearly, $N$ is almost solvable and
$G/N$ is abelian. But then, obviously, $G$ is almost solvable, which
is a contradiction.

Recall that by definition a group $G$ is Shunkov, if for any finite
subgroup $K$ of $G$, each subgroup of $N_{G}(K)/K$, generated by two
its conjugated elements of prime order, is finite. According to
Suchkova-Shunkov Theorem \cite{7}, which generalizes Shunkov's
Theorem \cite{1}, Shunkov groups with the minimal condition for
abelian subgroups are Chernikov. Also by \cite{8}, nonabelian
periodic Shunkov groups with the minimal condition for nonabelian
subgroups are Chernikov too.

Below we prove the following proposition based on Suchkova-Shunkov
Theorem \cite{7}.

\vskip2mm {\bf Theorem.} {\it Let $\mathfrak{Y}$ be a class of
groups in which any periodic subgroup is Shunkov, $\mathfrak{X}$ be
the minimal local class of groups containing $\mathfrak{Y}$ and
closed with respect to subgroups and ascending series. Then any
$\mathfrak{X}$-group satisfies the minimal condition for abelian
subgroups iff it is Chernikov. }

 \hspace{0.2cm}{\itshape\bfseries Proof. } {\it Sufficiency is obvious. }

 \hspace{0.2cm}{\it Necessity. } Let $\mathfrak{X_{0}}$ be the class consisting of all groups isomorphic
 to subgroups of $\mathfrak{Y}$-groups, and by induction for ordinals
 $\beta>0$: if there exists an ordinal $\alpha$ such that $\beta=\alpha+1$,
then $\mathfrak{X_{\beta}}$
 be the class of all groups, possessing a local system of subgroups that have a series with
 $\mathfrak{X}_{\alpha}$-factors, and if there is no such $\alpha$, then
 $\mathfrak{X}_{\beta}={\mathop{\bigcup}\limits_{\alpha<\beta}\mathfrak{X_{\alpha}}}$. It is easy to see that
 $\mathfrak{X}$ is the union of classes $\mathfrak{X}_{\beta}$.

 Suppose that the present theorem is not correct. Let $\gamma$ be the least one among all $\alpha$, for which  $\mathfrak{X}_{\alpha}$
 contains a non-Chernikov group $G_{\alpha}$ satisfying the minimal condition for abelian subgroups. In view of
 Suchkova-Shunkov Theorem \cite{7}, $\gamma>0$. It is easy to see that for some ordinal $\nu$, $\gamma=\nu+1$.
 Consequently, $G=G_{\gamma}$ has a local system $\mathcal{M}$ of subgroups possessing an ascending series
 with $\mathfrak{X_{\nu}}$-factors.

In view of Shunkov's Theorem \cite{1}, $G$ is not locally finite. So
some $H\in\mathcal{M}$ is not locally finite too. Let
$H_{0}=1\subset H_{1}\subset\ldots\subset H_{\theta}=H$ be an
ascending series of $H$ with $\mathfrak{X_{\nu}}$-factors and
$H_{\delta}$ be the union of all locally finite terms of the series.
Then, clearly, $H_{\delta}$ itself is locally finite. Since
$H_{\delta+1}$ is not locally finite, in consequence of
O.J.Schmidt's Theorem, $H_{\delta+1}/H_{\delta}$ is not locally
finite too. At the same time, $H_{\delta+1}/H_{\delta}$ is
non-Chernikov. The $\mathfrak{X_{\nu}}$-group
$H_{\delta+1}/H_{\delta}$ contains a non-Chernikov abelian subgroup
$A/H_{\delta}$. Since $H_{\delta+1}$ satisfies the minimal condition
for abelian subgroups, it is periodic. So $A/H_{\delta}$ is periodic
and, at the same time, locally finite. Then by O.J.Schmidt's
Theorem, $A$ is locally finite. Since $A$ satisfies the minimal
condition for abelian subgroups, it is Chernikov (Shunkov's Theorem
\cite{1}). But then $A/H_{\delta}$ must be Chernikov, which is a
contradiction. Theorem is proven.

\end{document}